\documentclass[12pt,vatola]{article}

\usepackage{graphicx}

\def\c{\centerline}

\def\re#1{\par\hangindent\parindent\indent\llap{#1\enspace}\ignorespaces}

\def\no{\noindent}

\begin{document}

\vskip 15mm

\c{\large\bf On Multi-Metric Spaces  }  \vskip 5mm

\c{Linfan Mao} \vskip 3mm \c{\scriptsize (Academy of Mathematics
and System Sciences, Chinese Academy of Sciences, Beijing 100080)}

\vskip 8mm
\begin{minipage}{130mm}
\no{\bf Abstract}: {\small A Smarandache multi-space is a union of
$n$ spaces $A_1,A_2,\cdots ,A_n$ with some additional conditions
holding. Combining Smarandache multi-spaces with classical metric
spaces, the conception of multi-metric space is introduced. Some
characteristics of a multi-metric space are obtained and Banach's
fixed-point theorem is generalized in this paper.}

\vskip 2mm \no{\bf Key words:} {\small  metric, multi-space,
multi-metric space, ideal subspace chain.}

 \vskip 2mm \no{{\bf
Classification:} AMS(2000) 51D99,51K05,51M05}
\end{minipage}

\vskip 8mm

{\bf $1$. Introduction}

\vskip 6mm

The notion of multi-spaces is introduced by Smarandache in $[6]$
under his idea of hybrid mathematics: {\it  combining different
fields into a unifying field}($[7]$), which is defined as follows.

\vskip 4mm

\no{\bf Definition $1.1$} \ {\it For any integer $i, 1\leq i\leq
n$ let $A_i$ be a set with ensemble of law $L_i$, and the
intersection of $k$ sets $A_{i_1},A_{i_2},\cdots , A_{i_k}$ of
them constrains the law $I(A_{i_1},A_{i_2},\cdots , A_{i_k})$.
Then the union of $A_i$, $1\leq i\leq n$

$$\widetilde{A} \ = \ \bigcup\limits_{i=1}^n A_i$$

\no is called a multi-space.}

\vskip 3mm

As we known, a set $M$ associative a function $\rho: M\times
M\rightarrow R^+=\{x \ | \ x\in R, x\geq 0\}$ is called a {\it
metric space} if for $\forall x,y,z\in M$, the following
conditions for the metric function $\rho$ hold:

\vskip 3mm

($1$)({\it definiteness}) \  $\rho (x,y)=0$ if and only if $x=y$;

($ii$)({\it symmetry}) \ $\rho (x,y)=\rho (y,x)$;

($iii$)({\it triangle inequality}) \ $\rho (x,y)+\rho (y,z)\geq
\rho (x,z).$

\vskip 2mm

By combining Smarandache multi-spaces with classical metric
spaces, a new kind of space called multi-ring space is found,
which is defined in the following.

\vskip 4mm

\no{\bf Definition $1.2$} \ {\it A multi-metric space is a union
$\widetilde{M}=\bigcup\limits_{i=1}^m M_i$ such that each $M_i$ is
a space with metric $\rho_i$ for $\forall i, 1\leq i\leq m$.}

When we say {\it a multi-metric space
$\widetilde{M}=\bigcup\limits_{i=1}^m M_i$}, it means that a
multi-metric space with metrics $\rho_1,\rho_2, \cdots ,\rho_m$
such that $(M_i,\rho_i)$ is a metric space for any integer $i,
1\leq i\leq m$. For a multi-metric space
$\widetilde{M}=\bigcup\limits_{i=1}^m M_i$, $x\in\widetilde{M}$
and a positive number $R$, a {\it $R$-disk} $B(x,R)$ in
$\widetilde{M}$ is defined by

$$B(x,R)=\{ \ y \ | \ {\rm there \ exists \ an \
integer} \ k, 1\leq k\leq m \ {\rm such \ that} \ \rho_k(y,x) \ <
R, y\in\widetilde{M}\}$$

The subject of this paper is to find some characteristics of a
multi-metric space. For terminology and notation not defined here
can be seen in $[1]-[2],[4]$ for terminologies in metric space and
in $[3],[5]-[9]$ for multi-spaces and logics.

\vskip 8mm

{\bf $2.$ Characteristics of a multi-metric space}

\vskip 5mm

For metrics on a space, we have the following result.

\vskip 4mm

\no{\bf Theorem $2.1$} \ {\it Let $\rho_1,\rho_2,\cdots ,\rho_m$
be $m$ metrics on a space $M$ and $F$ a function on ${\bf E}^m$
such that the following conditions hold:}

($i$) \ {\it $F(x_1,x_2,\cdots ,x_m)\geq F(y_1,y_2,\cdots ,y_m)$
if for $\forall i, 1\leq i\leq m$, $x_i\geq y_i$};

($ii$) \ {\it $F(x_1,x_2,\cdots ,x_m)=0$ only if $x_1=x_2=\cdots
=x_m=0$};

($iii$) \ {\it For two $m$-tuples $(x_1,x_2,\cdots ,x_m)$ and
$(y_1,y_2,\cdots ,y_m)$},

$$F(x_1,x_2,\cdots ,x_m)+F(y_1,y_2,\cdots ,y_m)\geq
F(x_1+y_1,x_2+y_2,\cdots ,x_m+y_m).$$

\no {\it Then $F(\rho_1,\rho_2,\cdots ,\rho_m)$ is also a metric
on $M$. }

\vskip 3mm

{\it Proof} \ We only need to prove that $F(\rho_1,\rho_2,\cdots
,\rho_m)$ satisfies the metric conditions for $\forall x,y,z\in
M$.

By ($ii$), $F(\rho_1(x,y),\rho_2(x,y),\cdots ,\rho_m(x,y))=0$ only
if for any integer $i$, $\rho_i(x,y)=0$. Since $\rho_i$ is a
metric on $M$, we know that $x=y$.

For any integer $i, 1\leq i\leq m$, since $\rho_i$ is a metric on
$M$, we know that $\rho_i(x,y)=\rho_i(y,x)$. Whence,

$$F(\rho_1(x,y),\rho_2(x,y),\cdots ,\rho_m(x,y))=
F(\rho_1(y,x),\rho_2(y,x),\cdots ,\rho_m(y,x)).$$

Now by ($i$) and ($iii$), we get that

\begin{eqnarray*}
& \ & F(\rho_1(x,y),\rho_2(x,y),\cdots
,\rho_m(x,y))+F(\rho_1(y,z),\rho_2(y,z),\cdots ,\rho_m(y,z))\\
& \ &\geq
F(\rho_1(x,y)+\rho_1(y,z),\rho_2(x,y)+\rho_2(y,z),\cdots,\rho_m(x,y)+\rho_m(y,z))\\
& \ &\geq F(\rho_1(x,z),\rho_2(x,z),\cdots ,\rho_m(x,z)).
\end{eqnarray*}

Therefore, $F(\rho_1,\rho_2,\cdots ,\rho_m)$ is a metric on $M$.
\quad\quad $\natural$

\vskip 4mm

\no{\bf Corollary $2.1$} \ {If $\rho_1,\rho_2,\cdots ,\rho_m$ are
$m$ metrics on a space $M$, then $\rho_1+\rho_2+\cdots +\rho_m$
and $\frac{\rho_1}{1+\rho_1}+\frac{\rho_2}{1+\rho_2}+\cdots
+\frac{\rho_m}{1+\rho_m}$ are also metrics on $M$.}

\vskip 3mm

A sequence $\{x_n\}$ in a multi-metric space
$\widetilde{M}=\bigcup\limits_{i=1}^m M_i$ is said {\it
convergence to a point} $x, x\in\widetilde{M}$ if for any positive
number $\epsilon > 0$, there exist numbers $N$ and $i, 1\leq i\leq
m$ such that if $n\geq N$ then

$$\rho_i(x_n,x) \ < \ \epsilon.$$

\no For $\{x_n\}$ convergence to a point $x, x\in\widetilde{M}$,
we denote it by $\lim\limits_{n}x_n = x$.

We have a characteristic for convergent sequences in a
multi-metric space.

\vskip 4mm

\no{\bf Theorem $2.2$} \ {\it A sequence $\{x_n\}$ in a
multi-metric space $\widetilde{M}=\bigcup\limits_{i=1}^m M_i$ is
convergent if and only if there exist integers $N$ and $k, 1\leq
k\leq m$, such that the subsequence $\{x_n| n\geq N\}$ is a
convergent sequence in $(M_k,\rho_k)$.}

\vskip 3mm

{\it Proof} \ If there exist integers $N$ and $k,1\leq k\leq m$,
such that $\{x_n | n\geq N\}$ is a convergent subsequence in
$(M_k,\rho_k)$, then for any positive number $\epsilon  > 0$, by
definition there exists a positive integer $P$ and a point $x,
x\in M_k$ such that

$$\rho_k(x_n,x) \ < \ \epsilon$$

\no if $n\geq max\{ N, \ P \}$.

Now if $\{x_n\}$ is a convergent sequence in the multi-space
$\widetilde{M}$, by definition for any positive number $\epsilon
> 0$,there exist a point $x, x\in\widetilde{M}$ and natural numbers
$N(\epsilon )$ and $k, 1\leq k\leq m$, such that if $n\geq
N(\epsilon)$, then

$$\rho_k(x_n,x) \ < \ \epsilon,$$

\no that is, $\{ x_n | n\geq N(\epsilon) \}\subset M_k$ and $\{
x_n | n\geq N(\epsilon) \}$ is a convergent sequence in
$(M_k,\rho_k)$.\quad\quad $\natural$

\vskip 4mm

\no{\bf Theorem $2.3$} \ {\it Let
$\widetilde{M}=\bigcup\limits_{i=1}^m M_i$ be a multi-metric
space. For two sequences $\{x_n\}, \{y_n\}$ in $\widetilde{M}$, if
$\lim\limits_{n}x_n = x_0$, $\lim\limits_{n}y_n = y_0$ and there
is an integer $p$ such that $x_0,y_0\in M_p$, then
$\lim\limits_n\rho_p(x_n,y_n)=\rho_p(x_0,y_0)$.}

\vskip 3mm

{\it Proof} \ According to Theorem $2.2$, there exist integers
$N_1$ and $N_2$ such that if $n\geq max\{N_1, N_2\}$, then
$x_n,y_n\in M_p$. Whence, we have that

$$\rho_p(x_n,y_n) \leq \rho_p(x_n,x_0)+\rho_p(x_0,y_0)+\rho_p(y_n,y_0)$$

\no{and}

$$\rho_p(x_0,y_0) \leq \rho_p(x_n,x_0)+\rho_p(x_n,y_n)+\rho_p(y_n,y_0).$$

\no Therefore,

$$|\rho_p(x_n,y_n)-\rho_p(x_0,y_0)| \ \leq \ \rho_p(x_n,x_0)+\rho_p(y_n,y_0).$$

For any positive number $\epsilon > 0$, since $\lim\limits_{n}x_n
= x_0$ and $\lim\limits_{n}y_n = y_0$, there exist numbers
$N_1(\epsilon), N_1(\epsilon)\geq N_1$ and $N_2(\epsilon),
N_2(\epsilon)\geq N_2$ such that $\rho_p(x_n,x_0) \ \leq \
\frac{\epsilon}{2}$ if $n\geq N_1(\epsilon)$ and $\rho_p(y_n,y_0)
\ \leq \ \frac{\epsilon}{2}$ if $n\geq N_2(\epsilon)$. Whence, if
$n\geq max\{N_1(\epsilon), N_2(\epsilon)\}$, then

$$|\rho_p(x_n,y_n)-\rho_p(x_0,y_0)| \ < \ \epsilon.\quad\quad \natural$$

Whether a convergent sequence can has more than one limit point?
The following result answers this question.

\vskip 4mm

\no{\bf Theorem $2.4$} \ {\it If $\{x_n\}$ is a convergent
sequence in a multi-metric space
$\widetilde{M}=\bigcup\limits_{i=1}^m M_i$, then $\{x_n\}$ has
only one limit point.}

\vskip 3mm

{\it Proof} \ According to Theorem $2.2$, there exist integers $N$
and $i, 1\leq i\leq m$ such that $x_n\in M_i$ if $n\geq N$. Now if

$$\lim\limits_nx_n = x_1 \ {\rm and} \ \lim\limits_nx_n = x_2,$$

\no and $n\geq N$, by definition,

$$0\leq\rho_i(x_1,x_2)\leq\rho_i(x_n,x_1)+\rho_i(x_n,x_2).$$

\no Whence, we get that $\rho_i(x_1,x_2)=0$. Therefore, $x_1=x_2$.
\quad\quad $\natural$

\vskip 4mm

\no{\bf Theorem $2.5$} \ {\it Any convergent sequence in a
multi-metric space is a bounded points set.}

\vskip 3mm

{\it Proof} \ According to Theorem $2.4$, we obtain this result
immediately. \quad\quad $\natural$

A sequence $\{x_n\}$ in a multi-metric space
$\widetilde{M}=\bigcup\limits_{i=1}^m M_i$ is called a {\it Cauchy
sequence} if for any positive number $\epsilon >0$, there exist
integers $N(\epsilon)$ and $s, 1\leq s\leq m$ such that for any
integers $m,n\geq N(\epsilon)$, $\rho_s(x_m,x_n) \ < \epsilon$.

\vskip 4mm

\no{\bf Theorem $2.6$} \ {\it A Cauchy sequence $\{x_n\}$ in a
multi-metric space $\widetilde{M}=\bigcup\limits_{i=1}^m M_i$ is
convergent if and only if for $\forall k, 1\leq k\leq m$,
$|\{x_n\}\bigcap M_k|$ is finite or infinite but $\{x_n\}\bigcap
M_k$ is convergent in $(M_k,\rho_k)$. }

\vskip 3mm

{\it Proof} \ The necessity of conditions is by Theorem $2.2$.

Now we prove the sufficiency. By definition, there exist integers
$s, 1\leq s\leq m$ and $N_1$ such that $x_n\in M_s$ if $n\geq
N_1$. Whence, if $|\{x_n\}\bigcap M_k|$ is infinite and
$\lim\limits_n\{x_n\}\bigcap M_k = x$, then there must be $k=s$.
Denoted by $\{x_n\}\bigcap M_k = \{x_{k1},x_{k2},\cdots ,x_{kn},
\cdots\}$.

For any positive number $\epsilon >0$, there exists an integer
$N_2, N_2\geq N_1$ such that $\rho_k(x_m,x_n) \ < \
\frac{\epsilon}{2}$ and $\rho_k(x_{kn},x) \ < \
\frac{\epsilon}{2}$ if $m,n\geq N_2$. According to Theorem $4.7$,
we get that

$$\rho_k(x_n,x)\leq \rho_k(x_n,x_{kn})+\rho_k(x_{kn},x)\ < \ \epsilon$$

\no if $n\geq N_2$. Whence, $\lim\limits_nx_n = x. \quad\quad
\natural$

A multi-metric space $\widetilde{M}$ is said {\it completed} if
every Cauchy sequence in this space is convergent. For a completed
multi-metric space, we obtain two important results similar to
metric space theory in classical mathematics.

\vskip 4mm

\no{\bf Theorem $2.7$} \ {\it Let
$\widetilde{M}=\bigcup\limits_{i=1}^m M_i$ be a completed
multi-metric space. For a $\epsilon$-disk sequence
$\{B(\epsilon_n,x_n)\}$, where $\epsilon_n > 0$ for $n=1,2,3,
\cdots$, the following conditions hold:}

($i$) \ {\it $B(\epsilon_1,x_1)\supset B(\epsilon_2,x_2)\supset
B(\epsilon_3,x_3)\supset\cdots\supset
B(\epsilon_n,x_n)\supset\cdots$};

($ii$) \ {\it $\lim\limits_{n\to+\infty}\epsilon_n=0$}.

\no {\it Then $\bigcap\limits_{n=1}^{+\infty}B(\epsilon_n,x_n)$
only has one point. }

\vskip 3mm

{\it Proof} \ First, we prove that the sequence $\{x_n\}$ is a
Cauchy sequence in $\widetilde{M}$. By the condition $(i)$, we
know that if $m\geq n$, then $x_m\in B(\epsilon_m, x_m)\subset
B(\epsilon_n, x_n)$. Whence, for $\forall i, 1\leq i\leq m$,
$\rho_i(x_m,x_n) \ < \epsilon_n$ if $x_m,x_n\in M_i$.

For any positive number $\epsilon$, since $\lim\limits_{n\to
+\infty}\epsilon_n=0$, there exists an integer $N(\epsilon)$ such
that if $n\geq N(\epsilon)$, then $\epsilon_n < \epsilon$.
Therefore, if $x_n\in M_l$, then $\lim\limits x_m=x_n$. Whence,
there exists an integer $N$ such that if $m\geq N$, then $x_m\in
M_l$ by Theorem $2.2$. Take integers $m,n\geq
max\{N,N(\epsilon)\}$. We know that

$$\rho_l(x_m,x_n) \ <  \ \epsilon_n \ < \ \epsilon.$$

\no So $\{x_n\}$ is a Cauchy sequence.

By the assumption, $\widetilde{M}$ is completed. We know that the
sequence $\{x_n\}$ is convergence to a point $x_0,
x_0\in\widetilde{M}$. By conditions ($i$) and ($ii$), we have that
$\rho_l(x_0,x_n) \ < \epsilon_n$ if we take $m\to +\infty$.
Whence, $x_0\in\bigcap\limits_{n=1}^{+\infty}B(\epsilon_n,x_n)$.

Now if there a point
$y\in\bigcap\limits_{n=1}^{+\infty}B(\epsilon_n,x_n)$, then there
must be $y\in M_l$. We get that

$$0\leq\rho_l(y,x_0) = \lim\limits_n\rho_l(y,x_n)\leq
\lim\limits_{n\to +\infty}\epsilon_n=0$$

\no by Theorem $2.3$. Therefore, $\rho_l(y,x_0)=0$. By definition
of a metric on a space, we get that $y=x_0$. \quad\quad $\natural$

Let  $\widetilde{M}_1$ and $\widetilde{M}_2$ be two multi-metric
spaces and $f: \widetilde{M}_1\rightarrow\widetilde{M}_2$ a
mapping, $x_0\in\widetilde{M}_1, f(x_0)=y_0$. For $\forall\epsilon
> 0$, if there exists a number $\delta$ such that for $forall x\in B(\delta ,
x_0)$, $f(x)=y\in B(\epsilon ,y_0)\subset\widetilde{M}_2$, i.e.,

$$f(B(\delta , x_0))\subset B(\epsilon ,y_0),$$

\no then we say $f$ is {\it continuous at point} $x_0$. If $f$ is
connected at every point of $\widetilde{M}_1$, then $f$ is said a
{\it continuous mapping} from $\widetilde{M}_1$ to
$\widetilde{M}_2$.

For a continuous mapping $f$ from $\widetilde{M}_1$ to
$\widetilde{M}_2$ and a convergent sequence $\{x_n\}$ in
$\widetilde{M}_1$, $\lim\limits_nx_n=x_0$, we can prove that

$$\lim\limits_nf(x_n)=f(x_0).$$

For a multi-metric space $\widetilde{M}=\bigcup\limits_{i=1}^m
M_i$ and a mapping $T: \widetilde{M}\rightarrow\widetilde{M}$, if
there is a point $x^*\in\widetilde{M}$ such that $Tx^* = x^*$,
then $x^*$ is called a {\it fixed point} of $T$. Denoted by
$^{\#}\Phi (T)$ the number of all fixed points of a mapping $T$ in
$\widetilde{M}$. If there are a constant $\alpha , 1 < \alpha < 1$
and integers $i, j, 1\leq i, j\leq m$ such that for $\forall
x,y\in M_i$, $Tx, Ty\in M_j$ and

$$\rho_j(Tx,Ty)\leq\alpha\rho_i(x,y),$$

\no then $T$ is called a {\it contraction} on $\widetilde{M}$.

\vskip 4mm

\no{\bf Theorem $2.8$} \ {\it Let
$\widetilde{M}=\bigcup\limits_{i=1}^m M_i$ be a completed
multi-metric space and $T$ a contraction on $\widetilde{M}$. Then}

$$1\leq ^{\#}\Phi (T)\leq m.$$

\vskip 3mm

{\it Proof} \ Choose arbitrary points $x_0, y_0\in M_1$ and define
recursively

$$ x_{n+1}=Tx_n, \ \ y_{n+1}=Tx_n$$

\no for $n=1,2,3,\cdots$. By definition, we know that for any
integer $n, n\geq 1$, there exists an integer $i, 1\leq i\leq m$
such that $x_n, y_n\in M_i$. Whence, we inductively get that

$$0\leq\rho_i(x_n,y_n)\leq\alpha^n\rho_1(x_0,y_0).$$

Notice that $0 < \alpha < 1$, we know that $\lim\limits_{n\to
+\infty}\alpha^n=0$. Therefore, there exists an integer $i_0$ such
that

$$\rho_{i_0}(\lim\limits_nx_n,\lim\limits_ny_n)=0.$$

Therefore, there exists an integer $N_1$ such that $x_n, y_n\in
M_{i_0}$ if $n\geq N_1$. Now if $n\geq N_1$, we have that

\begin{eqnarray*}
\rho_{i_0}(x_{n+1},x_n)&=& \rho_{i_0}(Tx_n,Tx_{n-1})\\
&\leq& \alpha
\rho_{i_0}(x_n,x_{n-1})=\alpha\rho_{i_0}(Tx_{n-1},Tx_{n-2})\\
&\leq& \alpha^2
\rho_{i_0}(x_{n-1},x_{n-2})\leq\cdots\leq\alpha^{n-N_1}\rho_{i_0}(x_{N_1+1},x_{N_1}).
\end{eqnarray*}

\no and generally, for $m\geq n\geq N_1$,

\begin{eqnarray*}
\rho_{i_0}(x_m,x_n) &\leq&
\rho_{i_0}(x_n,x_{n+1})+\rho_{i_0}(x_{n+1},x_{n+2})+\cdots
+\rho_{i_0}(x_{n-1},x_n)\\
&\leq& (\alpha^{m-1}+\alpha^{m-2}+\cdots
+\alpha^{n})\rho_{i_0}(x_{N_1+1},x_{N_1})\\
&\leq& \frac{\alpha^n}{1-\alpha}\rho_{i_0}(x_{N_1+1},x_{N_1})\to 0
(m,n\to +\infty)
\end{eqnarray*}

\no Therefore, $\{x_n\}$ is a Cauchy sequence in $\widetilde{M}$.
Similarly, we can prove $\{y_n\}$ is also a Cauchy sequence.

Because $\widetilde{M}$ is a completed multi-metric space, we have
that

$$\lim\limits_nx_n = \lim\limits_ny_n = z^*.$$

We prove $z^*$ is a fixed point of $T$ in $\widetilde{M}$. In
fact, by $\rho_{i_0}(\lim\limits_nx_n,\lim\limits_ny_n)=0$, there
exists an integer $N$ such that

$$x_n,y_n, Tx_n,Ty_n\in M_{i_0}$$

\no if $n\geq N+1$. Whence, we know that

\begin{eqnarray*}
0\leq\rho_{i_0}(z^*,Tz^*) &\leq& \rho_{i_0}(z^*,x_n)+
\rho_{i_0}(y_n,Tz^*)+\rho_{i_0}(x_n,y_n)\\
&\leq& \rho_{i_0}(z^*,x_n)+\alpha
\rho_{i_0}(y_{n-1},z^*)+\rho_{i_0}(x_n,y_n).
\end{eqnarray*}

\no Notice

$$\lim\limits_{n\to +\infty}\rho_{i_0}(z^*,x_n)=\lim\limits_{n\to
+\infty}\rho_{i_0}(y_{n-1},z^*)=\lim\limits_{n\to
+\infty}\rho_{i_0}(x_n,y_n)=0.$$

\no We get that $\rho_{i_0}(z^*,Tz^*)=0$, i.e., $Tz^*=z^*$.

For other chosen points $u_0,v_0\in M_1$, we can also define
recursively

$$u_{n+1} = Tu_n, \ \ v_{n+1} = Tv_n$$

\no and get the limit points $\lim\limits_nu_n = \lim\limits_nv_n
= w^*\in M_{i_0}, Tu^*\in M_{i_0}$. Since

$$\rho_{i_0}(z^*,u^*)=\rho_{i_0}(Tz^*,Tu^*)\leq\alpha\rho_{i_0}(z^*,u^*)$$

\no and $0 < \alpha < 1,$ there must be $z^*=u^*$.

Similar consider the points in $M_i, 2\leq i\leq m$, we get that

$$1\leq ^{\#}\Phi (T)\leq m. \quad\quad\quad \natural$$

\vskip 4mm

\no{\bf Corollary $2.2$}(Banach) \ {\it Let $M$ be a metric space
and $T$ a contraction on $M$. Then $T$ has just one fixed point.}

\vskip 8mm

{\bf $3.$ Open problems for a multi-metric space}

\vskip 5mm

On a classical notion, only one metric maybe considered in a space
to ensure the same on all the times and on all the situations.
Essentially, this notion is based on an assumption that spaces are
homogeneous. In fact, it is not true in general.

Multi-Metric spaces can be used to simplify or beautify
geometrical figures and algebraic equations. One example is shown
in Fig.$1$, in where the left elliptic curve is transformed to the
right circle by changing the metric along $x, y$-axes and an
elliptic equation

$$\frac{x^2}{a^2}+\frac{y^2}{b^2}=1$$

\no to equation

$$x^2+y^2=r^2$$

\no of a circle of radius $r$.

\vskip 2mm

\includegraphics[bb=75 5 400 190]{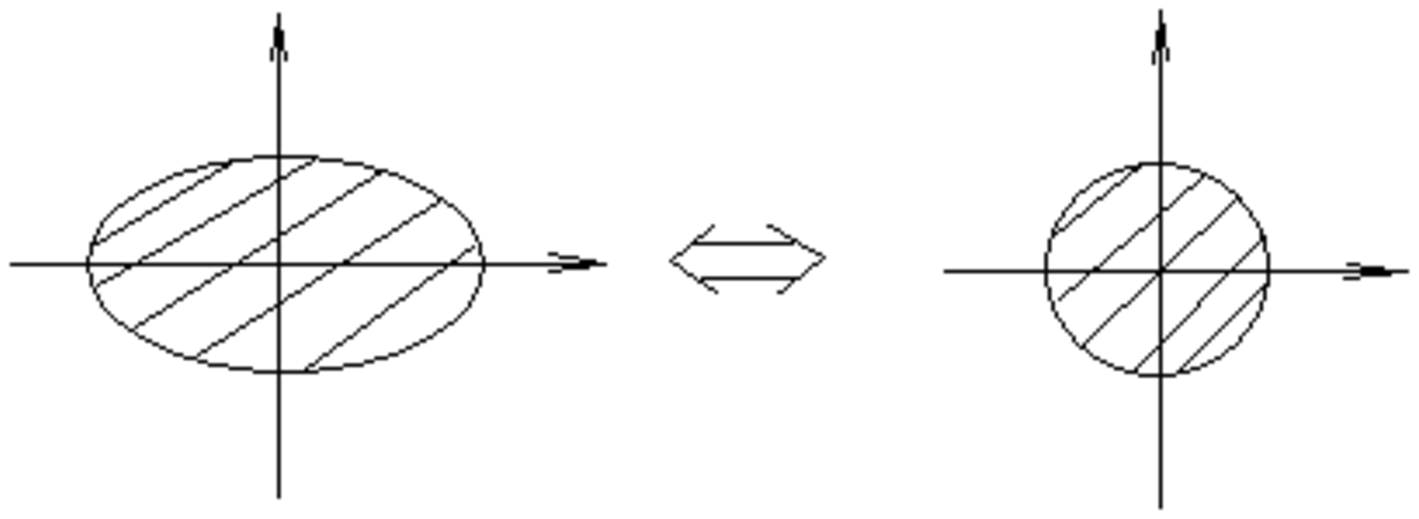}

\vskip 2mm

\c{\bf Fig.$1$}\vskip 3mm

Generally, in a multi-metric space we can simplify a polynomial
similar to the approach used in projective geometry. {\it Whether
this approach can be contributed to mathematics with metrics?}

\no{\bf Problem $3.1$} \ {\it Choose suitable metrics to simplify
the equations of surfaces and curves in ${\bf E}^3$.}

\vskip 3mm

\no{\bf Problem $3.2$} \ {\it Choose suitable metrics to simplify
the knot problem. Whether can it be used for classifying
$3$-dimensional manifolds?}

\vskip 3mm

\no{\bf Problem $3.3$} \ {\it Construct multi-metric spaces or
non-linear spaces by Banach spaces. Simplify equations or problems
to linear problems.}

\vskip 8mm

{\bf References}

\vskip 5mm

\re{[1]}R. Abraham, J.E.Marsden and T.Ratiu, {\it Manifolds,
Tensor Analysis and Applications}, Addison-Wesley Publishing
Company, Inc., 1983.

\re{[2]}S.S.Chern and W.H.Chern, {\it Lectures in Differential
Geometry}, Peking University Press, 2001.

\re{[3]}Daniel Deleanu, {\it A Dictionary of Smarandache
Mathematics}, Buxton University Press, London \& New York,2004.

\re{[4]}J.M.Lee, {\it Riemannian Manifolds}, Springer-Verlag New
York, Inc.,1997.

\re{[5]}L.F.Mao, {\it Automorphism Groups of Maps, Surfaces and
Smarandache Geometries}, American Research Press, 2005.

\re{[6]} F.Smarandache, Mixed noneuclidean geometries, {\it eprint
arXiv: math/0010119}, 10/2000.

\re{[7]}F.Smarandache, {\it A Unifying Field in Logics.
Neutrosopy: Neturosophic Probability, Set, and Logic}, American
research Press, Rehoboth, 1999.

\re{[8]}F.Smarandache, Neutrosophy, a new Branch of Philosophy,
{\it Multi-Valued Logic}, Vol.8, No.3(2002)(special issue on
Neutrosophy and Neutrosophic Logic), 297-384.

\re{[9]}F.Smarandache, A Unifying Field in Logic: Neutrosophic
Field, {\it Multi-Valued Logic}, Vol.8, No.3(2002)(special issue
on Neutrosophy and Neutrosophic Logic), 385-438.

\end{document}